\documentclass[letterpaper, 10 pt, conference]{ieeeconf}

\IEEEoverridecommandlockouts
\usepackage[utf8]{inputenc}
\usepackage[T1]{fontenc}
\usepackage[english]{babel}

\usepackage[rgb]{xcolor}
\usepackage{microtype}                      
\usepackage[xlatin, abbreviations]{foreign} 

\usepackage{amsmath, amssymb}               
\usepackage{amsthm}
\usepackage{mathtools}                      
\usepackage[casesstyle, subnum]{cases}      
\usepackage{bm}                             

\usepackage{enumerate}                      
\usepackage{booktabs}                       
\usepackage{siunitx}

\makeatletter
\let\NAT@parse\undefined
\makeatother

\usepackage{csquotes}                       
\usepackage{hyperref}
\usepackage{cite}                           

\usepackage{tikz, tikz-3dplot}
\usepackage{pgfplots}                       
\usepackage{pgfplotstable}
\usepackage[caption=false,font=footnotesize]{subfig}

\usetikzlibrary{arrows.meta, calc, patterns, positioning}

\usepgfplotslibrary{colorbrewer}            
\pgfplotsset{compat=1.18}

\newtheorem{theorem}{Theorem}
\newtheorem{corollary}[theorem]{Corollary}
\newtheorem{lemma}[theorem]{Lemma}
\newtheorem{proposition}[theorem]{Proposition}

\theoremstyle{definition}
\newtheorem{assumption}{Assumption}
\newtheorem{definition}{Definition}

\theoremstyle{remark}
\newtheorem*{remark}{Remark}

\NewDocumentCommand{\real}{}{\mathbb{R}}

\DeclareMathOperator{\binary}{bin}

\DeclareMathOperator{\relint}{relint}

\DeclareMathOperator{\conv}{conv}

\ProvideDocumentCommand{\given}{}{}
\NewDocumentCommand{\renewgiven}{O{\delimsize}}{\renewcommand\given{\nonscript\:#1\vert\allowbreak\nonscript\:\mathopen{}}}

\DeclarePairedDelimiterX\card[1]\lvert\rvert{#1}
\DeclarePairedDelimiterX\norm[1]\lVert\rVert{#1}
\DeclarePairedDelimiterX\scprod[2]\langle\rangle{#1, #2}
\DeclarePairedDelimiterXPP\diag[1]{\operatorname{diag}}(){}{#1}

\DeclarePairedDelimiterX\set[1]\{\}{\renewgiven #1}
\DeclarePairedDelimiterXPP\prob[1]{\operatorname{\mathbb{P}}}(){}{\renewgiven #1}
\DeclarePairedDelimiterXPP\expect[1]{\operatorname{\mathbb{E}}}[]{}{\renewgiven #1}
\DeclarePairedDelimiterXPP\variance[1]{\operatorname{Var}}[]{}{\renewgiven #1}

\makeatletter
\renewcommand*\env@matrix[1][*\c@MaxMatrixCols c]{%
	\hskip -\arraycolsep
	\let\@ifnextchar\new@ifnextchar
	\array{#1}}
\newcommand*{\matline}[1]{\cmidrule(lr){1-#1}}
\makeatother

\NewDocumentCommand{\arxivTF}{m m}{#1}

\usepackage{fancyhdr}

\fancypagestyle{firstpage}{%
    \fancyfoot[C]{\small
        © 2024 IEEE.
        Personal use of this material is permitted.
        Permission from IEEE must be obtained for all other uses, in any current or future media, including reprinting/republishing this material for advertising or promotional purposes, creating new collective works, for resale or redistribution to servers or lists, or reuse of any copyrighted component of this work in other works.
    }%
}

\title{\LARGE \bf
Robust Stability for Multiagent Systems with \\ Spatio-Temporally Correlated Packet Loss
}

\author{Christian Hespe and Adwait Datar and Herbert Werner%
\thanks{C. Hespe, A. Datar, and H. Werner are with the Hamburg University of Technology, Institute of Control Systems, Eißendorfer Straße 40, 21073 Hamburg, Germany, {\tt\small \{christian.\allowbreak hespe, adwait.\allowbreak datar, h.werner\}@tuhh.de}}%
}

\begin{document}

\maketitle
\thispagestyle{firstpage}

\begin{abstract}
    A problem with considering correlations in the analysis of multiagent system with stochastic packet loss is that they induce dependencies between agents that are otherwise decoupled, preventing the application of decomposition methods required for efficient evaluation.
To circumvent that issue, this paper is proposing an approach based on analysing sets of networks with independent communication links, only considering the correlations in an implicit fashion.
Combining ideas from the robust stabilization of Markov jump linear systems with recently proposed techniques for analysing packet loss in multiagent systems, we obtain a linear matrix inequality based stability condition which is independent of the number of agents.
The main result is that the set of stabilized probability distributions has non-empty interior such that small correlations cannot lead to instability, even though only distributions of independent links were analysed.
Moreover, two examples are provided to demonstrate the applicability of the results to practically relevant scenarios.

\end{abstract}

\pagestyle{empty}
\section{Introduction}\label{sec:intro}
Interconnected systems and multiagent systems (MASs) in particular have spurred a lot of interest in the last two decades since they are uniquely suited to solve many large-scale control problems.
Whereas centralized solutions fail to scale beyond a certain size, the distributed or decentralized nature of interconnected systems allows for an efficient scaling of the implementation \cite{MesbaE2010}.

An important aspect in the implementation of interconnected systems is how the information exchange between subsystems is handled.
Especially for MASs, the choice often falls on wireless communication networks.
However, even though wireless networks are inherently unreliable \cite{SchenSFPS2007}, the communication aspect is neglected in many works on MASs.
For example, this is the case in the decomposable systems framework proposed by Massioni and Verhaegen in \cite{MassiV2009}, which allows not only for scalable implementation but also scalable system analysis and controller synthesis.
Only recently, \cite{HespeSDWT2023} extended the framework to cover Markov jump linear system (MJLS) and thus consider stochastic packet loss and its detrimental effects.

Amongst the works that \emph{do} consider the stochastic effects introduced by wireless networking in MASs, one feature is a common necessity:
The need for simplifying assumptions on the kind of packet loss the network is subject to.
This is caused by the need for some form of decoupling in the analysis in order to handle the scale of MASs, while the most general stochastic communication models with arbitrary correlations would introduce coupling between all agents.
In the simplest case, studied amongst others in \cite{ZhangTHK2017, XuMX2020}, it is assumed that all communication links fail or succeed simultaneously, resulting in the need for only a single random variable.
For most practical systems, this assumption is too restrictive to be applicable, therefore other works consider a set of \emph{independent} stochastic variables, one for each edge.
Often, it is then assumed that the loss probability is homogeneous across all communication links, \eg in \cite{MesbaE2010, HespeSDWT2023, PatteB2010} and \cite{WuS2012}, again limiting which system the results can be applied to.
Approaches that allow for heterogeneous probabilities are less frequent.
For example, \cite{ZhangT2012} and \cite{HespeW2023} handle such cases by introducing homogeneous bounds on the probabilities, while \cite{MaXL2020,GordoVP2023} restrict the analysis to directed tree graphs.

The works above have in common that they rely on either independence or full correlation of the packet loss on individual links.
Two approaches that explicitly consider the probability distribution of the whole communication network are \cite{XuMX2020} and \cite{ChenKZC2021}.
However, both rely on constructing the transition probability matrix (TPM) in their analysis, which is prohibitively expensive for large MASs, since the TPM is in general a $2^m \times 2^m$ matrix for a MAS with $m$ communication links.
In contrast, the current paper proposes to analyse a collection of probability distribution with independent communication links, neglecting the correlations in the network at first.
For this kind of packet loss distribution, we can combine the analysis results for MJLSs with uncertain TPM presented in \cite{OliveVVP2009}, which by themselves are unsuitable for systems with many modes, with the recently proposed approach of \cite{HespeW2023} to obtain scalable analysis conditions for a subset of MJLS called \emph{decomposable MJLS} that is suitable for many MASs.
The key result is shown afterwards, by proving that the considered kind of uncorrelated packet loss distributions are in the interior of the uncertainty set relative to \emph{all possible} distributions including those with correlations and, as such, there cannot exist small correlations that induce instability.
In this way, we can take advantage of the excellent scalability of the results presented in \cite{HespeW2023} even though the underlying distribution does not necessarily feature independent links.

After the introduction, Section~\ref{sec:mas} defines the system and packet loss model considered in this paper.
The main robust stability test is derived in Section~\ref{sec:uncertain}, while Section~\ref{sec:simplex} is proving geometric properties of the involved uncertainty set.
To demonstrate the applicability of the proposed results, Section~\ref{sec:example} is discussing two examples, before Section~\ref{sec:conclusions} concludes the paper with closing remarks.

\subsection{Contribution}\label{sec:intro_contribution}
The contribution of this paper is a linear matrix inequality (LMI) based sufficient condition for distributionally robust stability of MASs subject to spatio-temporal correlations between communication links.
It consists of two main components:
\begin{enumerate}
    \item The stability test in Theorem~\ref{thm:bernoulli_vertices} in terms of a set of uncertain TPMs, and
    \item a characterization of the uncertainty set in Theorem~\ref{thm:uncertainty_richness}, proving that it contains the considered kind of independent distributions in its relative interior.
\end{enumerate}
An important feature is the independence of its computational complexity from the number of agents, enabling analysis of arbitrarily large MASs.

\subsection{Definitions \& Notation}\label{sec:intro_notation}
We use $I_N$ and $\mathbf{1}_N$ to denote the $N \times N$ identity matrix and $N$-dimensional vector of ones, respectively.
$A \otimes B$ is the Kronecker product of $A$ and $B$.
$M \succ 0$ and $M \prec 0$ mean that the matrix $M$ is positive or negative definite and $*$ is used to indicate entries required for symmetry.
A matrix is row-stochastic if its entries are non-negative and its rows sum up to one.
Furthermore, for a set $\mathcal{M}$, $\conv\mathcal{M}$ denotes its convex hull and $\relint\mathcal{M}$ its relative interior.

A graph $\mathcal{G} \coloneqq (\mathcal{V}, \mathcal{E})$ is composed of the vertex set $\mathcal{V} \coloneqq \set{1, 2, \dotsc, N}$ and the edge set $\mathcal{E} \subset \mathcal{V} \times \mathcal{V}$.
In this paper, we consider only undirected graphs such that an edge $e^{ij} \equiv e^{ji} \coloneqq \set{i,j}$ is an unordered pair of vertices and is seen as a bidirectional connection between vertices $i$ and $j$.
Furthermore, we assume no self loops exist, \ie $e^{ii} \notin \mathcal{E}$.
Finally, the Laplacian is defined as $L(\mathcal{G}) = [l_{ij}(\mathcal{G})]$ with
\begin{equation*}
    l_{ij}(\mathcal{G}) \coloneqq \begin{cases}
        -1   & if \(i \neq j\) and \(e_{ij} \in \mathcal{E}\), \\
         0   & if \(i \neq j\) and \(e_{ij} \notin \mathcal{E}\), \\
        -\textstyle{\sum_{l \neq i}} l_{il}(\mathcal{G}) & if \(i = j\).
    \end{cases}
\end{equation*}

\section{Multiagent Systems with Packet Loss}\label{sec:mas}
\subsection{System \& Packet Loss Model}\label{sec:mas_model}
Consider a MAS with $N \geq 2$ agents that are exchanging information over an unreliable communication network, which we are modelling using tools from graph theory.
The nominal interconnection topology, \ie if no information is lost, is captured by the graph $\mathcal{G}^0 \coloneqq (\mathcal{V}, \mathcal{E}^0)$ with one-to-one correspondence between vertices and agents.
We define $m \coloneqq \card{\mathcal{E}^0}$ and furthermore introduce a stochastic process $\{\theta_k(e^{ij})\}$ for every edge $e^{ij} \in \mathcal{E}^0$ that takes values in $\set{0,1}$ and governs the reception and loss of information.
At step $k$, $\theta_k(e^{ij}) = 1$ means a packet is successfully transmitted over $e^{ij}$ and correspondingly data is lost for $\theta_k(e^{ij}) = 0$.

By defining a function $\nu: \mathcal{E}^0 \to \set{1, 2, \dotsc, m}$ that associates every $e^{ij} \in \mathcal{E}^0$ with a unique integer, we combine all $\{\theta_k\}$ into a single stochastic process
\begin{equation}\label{eq:sigma}
    \sigma_k = 1 + \sum_{e^{ij} \in \mathcal{E}^0} \theta_k(e^{ij}) 2^{m-\nu(e^{ij})}
\end{equation}
with $\sigma_k \in \mathcal{K} \coloneqq \set[\big]{1, 2, \dotsc, 2^m}$.
For modelling purposes, we then make the following assumption:

\begin{assumption}\label{ass:markov}
    The stochastic process $\{\sigma_k\}$ is a homogeneous Markov chain, \ie there exist $t_{ij} \geq 0$ such that
    \begin{equation*}
        \prob{\sigma_{k+1} = j \given \sigma_k = i} = t_{ij}
    \end{equation*}
    for all $k \geq 0$ and $i,j \in \mathcal{K}$, where $\sum_{l \in \mathcal{K}} t_{il} = 1$ for all $i \in \mathcal{K}$.
\end{assumption}

In addition to the nominal graph $\mathcal{G}^0$, the packet loss is inducing graphs $\mathcal{G}_i = (\mathcal{V}, \mathcal{E}_i)$ with $i \in \mathcal{K}$, where $\mathcal{E}_i \subseteq \mathcal{E}^0$ is the subset of edges that successfully transmit information in mode $i$.
Combined with the Markov chain $\{\sigma_k\}$, these graphs are determining the dynamics of the MAS at a particular time step $k$.
This behaviour is described by the MJLS
\begin{equation}\label{eq:mjls}
    x_{k+1} = A_{\sigma_k} x_k,
\end{equation}
where $x_k \in \real^{N n_x}$ is the dynamic state and the system matrix is switched amongst $\set{A_i \in \real^{N n_x \times N n_x} \given i \in \mathcal{K}}$ (\cf \cite{CostaMF2005}).
More specifically, this paper is focused on systems that fit into the decomposable systems framework proposed by Massioni and Verhaegen in \cite{MassiV2009} and extended to MJLS in \cite{HespeSDWT2023}.
Therefore, we assume that the system matrix is structured as
\begin{equation}\label{eq:decomposable}
    A_i = I_N \otimes A^d + L(\mathcal{G}_i) \otimes A^c,
\end{equation}
where $A^c$ and $A^d$ are the coupled and decoupled component, respectively, and the Laplacian takes the role of the pattern matrix.
In the following, we will be using the shorthand notation $L^0 \coloneqq L(\mathcal{G}^0)$ and $L_i \coloneqq L(\mathcal{G}_i)$ if it is clear from context which graph the Laplacian corresponds to.

\subsection{Independent Packet Loss Distributions}\label{sec:mas_independence}
In contrast to most existing work on MASs with stochastic packet loss, we make no \apriori assumption on the spatial or temporal independence of $\{\theta_k(e^{ij})\}$ except for jointly forming a homogeneous Markov chain.
However, analysing the system with a packet loss model this general is computationally intractable even for moderately sized MAS \cite{HespeSDWT2023}.
Instead, we propose to use the approach recently introduced in \cite{HespeW2023} to analyse the MAS with a simplified, spatially independent packet loss model and show that these results are distributionally robust, including against spatio-temporal correlation.

\begin{figure}
    \centering
    \begin{tikzpicture}[
        >=latex,
        agent/.style={
            draw, circle, fill=gray!30, inner sep=2pt
        },
    ]
    \draw (0,0) node[agent](v1){$1$};
    \draw node[agent, right=of v1](v2){$2$};
    \draw node[agent, right=of v2](v3){$3$};
    
    \draw[<->] (v1) -- node[above]{$e^{12}$}(v2);
    \draw[<->] (v2) -- node[above]{$e^{23}$}(v3);
\end{tikzpicture}
    \caption{Exemplary multiagent system with three agents and two links}
    \label{fig:three_agents}
\end{figure}
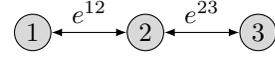

As an example, consider the MAS of three agents shown in Fig.~\ref{fig:three_agents}.
By the packet loss model introduced in the previous subsection, we consider the joint probability distribution of $\theta_k(e^{12})$ and $\theta_k(e^{23})$ which forms the three-dimensional simplex (one dimension is lost due to the constraint that the joint probabilities have to sum to one) .
If, instead, we would assume independence of $\theta_k(e^{12})$ and $\theta_k(e^{23})$, we obtain
\begin{align*}
    \begin{bmatrix}
        \prob{\theta_k(e^{12}) = 0, \theta_k(e^{23}) = 0} \\
        \prob{\theta_k(e^{12}) = 0, \theta_k(e^{23}) = 1} \\
        \prob{\theta_k(e^{12}) = 1, \theta_k(e^{23}) = 0} \\
        \prob{\theta_k(e^{12}) = 1, \theta_k(e^{23}) = 1}
    \end{bmatrix} = \begin{bmatrix}
        (1 - p^{12}) (1 - p^{23}) \\
        (1 - p^{12}) p^{23} \\
        p^{12} (1 - p^{23}) \\
        p^{12} p^{23}
    \end{bmatrix}
\end{align*}
with $p^{ij} \coloneqq \prob[\big]{\theta_k(e^{ij}) = 1}$ and are thus left with only two degrees of freedom.
This independence model thus covers a two-dimensional surface embedded in the three-dimensional simplex \cite[Fig.~1]{Montufar2013}.
However, we argue that by robustly stabilizing a sufficiently large number of distributions with independent $\theta_k(e^{12})$ and $\theta_k(e^{23})$, we obtain a stability guarantee that extends to distributions featuring spatial correlation.
This statement is made precise in Section~\ref{sec:simplex}.

\section{Uncertain Transition Probabilities}\label{sec:uncertain}
\subsection{Distributionally Robust MJLS Analysis}\label{sec:uncertain_distribution}
Suppose now that the transition probabilities are not known precisely, but that the TPM $\Gamma = [t_{ij}]$ belongs to the set
\begin{equation}\label{eq:tpm_set}
    \bm{\Gamma} \coloneqq \set*{\Gamma \in \real^{m \times m} \given \Gamma = \sum_{r = 1}^{n_\Gamma} \alpha_r \Gamma^{(r)}, \alpha \in \mathcal{S}_{n_\Gamma}},
\end{equation}
where $\mathcal{S}_n$ is the unit simplex given by
\begin{equation}\label{eq:simplex}
    \mathcal{S}_n \coloneqq \set[\big]{\alpha \in \real^n \given \textstyle\sum_{i=1}^n \alpha_i = 1, \alpha_i \geq 0}
\end{equation}
and $\Gamma^{(r)} = \bigl[t_{ij}^{(r)}\bigr]$ are known row-stochastic matrices.
The TPMs in $\bm{\Gamma}$ are thus parametrized by $\alpha \in \mathcal{S}_{n_\Gamma}$ and we write $\Gamma(\alpha)$ -- and accordingly $t_{ij}(\alpha)$ -- to refer to the TPM corresponding to a specific value of $\alpha$.

Based on the approach in \cite{OliveVVP2009}, we want to study under which conditions the MJLS is stable for all admissible TPMs.
First, we define stability for the MJLS:
\begin{definition}[Robust Mean-Square Stability]
    The MJLS~\eqref{eq:mjls} is said to be robustly mean-square stable if
    \begin{align*}
        \lim_{k \to \infty} \expect{x_k} &= 0 & &\text{and} & \lim_{k \to \infty} \expect[\big]{\norm{x_k}^2} &= 0
    \end{align*}
    for all $x_0 \in \real^n, \sigma_0 \in \mathcal{K}$ and $\Gamma \in \bm{\Gamma}$.
\end{definition}

The results of \cite{OliveVVP2009} already include an LMI based test for robust mean-square stability of the MJLS~\eqref{eq:mjls}.
However, for our purpose an equivalent form of their result is more suitable, which we state next:
\begin{theorem}\label{thm:mjls_stability}
    The MJLS~\eqref{eq:mjls} is robustly mean-square stable if and only if there exist functions $X_i(\alpha) \succ 0$ such that
    \begin{equation}\label{eq:mjls_stability_lmi}
        \sum_{j \in \mathcal{K}} t_{ij}(\alpha) A_j^\top X_j(\alpha) A_j - X_i(\alpha) \prec 0
    \end{equation}
    holds for all $i \in \mathcal{K}$ and $\alpha \in \mathcal{S}_{n_\Gamma}$.
\end{theorem}
\begin{proof}
    \arxivTF{%
        The original conditions from \cite[Lemma 3]{OliveVVP2009} demand that there exist functions $P_i(\alpha) \succ 0$ such that 
\begin{equation}\label{eq:mjls_stability_equiv_lmi}
        A_i^\top \left(\sum\nolimits_{j \in \mathcal{K}} t_{ij}(\alpha) P_j(\alpha)\right) A_i - P_i(\alpha) \prec 0
\end{equation}
holds for all $i \in \mathcal{K}$ and $\alpha \in \mathcal{S}_{n_\Gamma}$.
To apply the procedure from \cite{FioraGDG2012}, we first assume there exist functions $P_i(\alpha) \succ 0$ solving \eqref{eq:mjls_stability_equiv_lmi} for all $i \in \mathcal{K}$ and define $X_i(\alpha) \coloneqq \sum_{j \in \mathcal{K}} t_{ij}(\alpha) P_j(\alpha)$.
Then, for each $j \in \mathcal{K}$, multiply \eqref{eq:mjls_stability_equiv_lmi} with $t_{ji}(\alpha)$ and sum over $i \in \mathcal{K}$ to obtain
\begin{align*}
    0 &\succ \sum\nolimits_{i \in \mathcal{K}} t_{ji}(\alpha) \left[A_i^\top \left(\sum\nolimits_{j \in \mathcal{K}} t_{ij}(\alpha) P_j(\alpha)\right) A_i - P_i(\alpha)\right] \\
    &= \sum\nolimits_{i \in \mathcal{K}} t_{ji}(\alpha) A_i^\top X_i(\alpha) A_i - X_j(\alpha),
\end{align*}
which is \eqref{eq:mjls_stability_lmi} with swapped $i,j$.

On the other hand, assume there are $X_i(\alpha) \succ 0$ that solve \eqref{eq:mjls_stability_lmi} and define $P_i(\alpha) \coloneqq A_i^\top X_i(\alpha) A_i + \varepsilon I \succ 0$ with $\varepsilon > 0$.
For sufficiently small $\varepsilon > 0$, we thus have
\begin{multline*}
    \sum\nolimits_{j \in \mathcal{K}} t_{ij}(\alpha) P_j(\alpha) - X_i(\alpha) \\
    = \sum\nolimits_{j \in \mathcal{K}} t_{ij}(\alpha) A_j^\top X_j(\alpha) A_j - X_i(\alpha) + \varepsilon I
    \prec 0
\end{multline*}
for all $i \in \mathcal{K}$.
Finally, it follows that
\begin{multline*}
    A_i^\top \left(\sum\nolimits_{j \in \mathcal{K}} t_{ij}(\alpha) P_j(\alpha)\right) A_i - P_i(\alpha) \\
    \preceq A_i^\top X_i(\alpha) A_i - P_i(\alpha) = -\varepsilon I \prec 0
\end{multline*}
for all $i \in \mathcal{K}$.
    }{%
        The equivalence to \cite[Lemma~3]{OliveVVP2009} is established using the procedure from \cite{FioraGDG2012}.
        See \cite{ExtendedVersion} for details.
    }
\end{proof}

The stability test in Theorem~\ref{thm:mjls_stability} is an infinite-dimensional feasibility problem as the variables are arbitrary matrix-valued functions and the constraints have to be satisfied for all $\alpha \in \mathcal{S}_{n_\Gamma}$.
The problem therefore has to be simplified to be numerically tractable, which we achieve by restricting the search to constant matrices.
\begin{corollary}\label{cor:mjls_stability_vertices}
    The MJLS~\eqref{eq:mjls} is robustly mean-square stable if there exist $X_i \succ 0$ such that
    \begin{equation}\label{eq:mjls_stability_vertices_lmi}
        \sum_{j \in \mathcal{K}} t_{ij}^{(r)} A_j^\top X_j A_j - X_i \prec 0
    \end{equation}
    holds for all $i \in \mathcal{K}$ and $r \in \set{1, \ldots, n_\Gamma}$.
\end{corollary}
\begin{proof}
    From \eqref{eq:tpm_set}, it follows that $t_{ij}(\alpha) = \sum_{r = 1}^{n_\Gamma} \alpha_r t_{ij}^{(r)}$.
Substituting into \eqref{eq:mjls_stability_lmi}, we obtain
\begin{multline*}
    \sum_{j \in \mathcal{K}} \Biggl(\sum_{r = 1}^{n_\Gamma} \alpha_r t_{ij}^{(r)}\Biggr) A_j^\top X_j A_j - X_i \\
    = \sum_{r = 1}^{n_\Gamma} \alpha_r \Biggl(\sum_{j \in \mathcal{K}} t_{ij}^{(r)} A_j^\top X_j A_j - X_i\Biggr)
    \prec 0,
\end{multline*}
where we have used that $X_i = \sum_{r=1}^{n_\Gamma} \alpha_r X_i$ and the inequality is implied by \eqref{eq:mjls_stability_vertices_lmi} for all $r \in \set{1, \ldots, n_\Gamma}$.

\end{proof}
\begin{remark}
    Restricting Corollary~\ref{cor:mjls_stability_vertices} to constant matrices is a particularly simple but conservative choice.
    A less conservative reformulation is obtained in \cite{OliveVVP2009} by allowing $X_i(\alpha)$ to be homogeneous polynomials and applying Pólya relaxations of increasing order.
    However, these reformulations do not scale well to systems with large $n_\Gamma$, which will be the case for the MAS considered in this paper.
\end{remark}

\subsection{Spatially Independent Vertices}\label{sec:uncertain_vertices}
Corollary~\ref{cor:mjls_stability_vertices} shows that sufficient stability conditions can be obtained by testing only the vertices of $\bm{\Gamma}$.
Nonetheless, testing \eqref{eq:mjls_stability_lmi} requires enumerating $j \in \mathcal{K}$, which is intractable for MAS with more than a few agents because $\card{\mathcal{K}}$ grows exponentially with the number of edges.
We therefore take the approach of carefully choosing the vertices $\Gamma^{(r)}$ in a way such that they can be efficiently analysed.

The first simplification is to restrict the matrix variables to be identical and block-repeated, \ie $X_i = X_j = I_N \otimes \tilde{X}$.
Due to this change, we can take advantage of the mixed-product property of the Kronecker product \cite{Steeb1991} to obtain
\begin{equation}\label{eq:mjls_stability_decomposed_lmi}\begin{split}
    I_N& \otimes \bigl(A^{d\top} \tilde{X} A^d - \tilde{X}\bigr) \\
    &+ \Bigl(\sum\nolimits_{j \in \mathcal{K}} t_{ij}^{(r)} L_j^2 \Bigr) \otimes (A^{c\top} \tilde{X} A^c\bigr) \\
    &+ \Bigl(\sum\nolimits_{j \in \mathcal{K}} t_{ij}^{(r)} L_j \Bigr) \otimes \bigl(A^{d\top} \tilde{X} A^c + A^{c\top} \tilde{X} A^d\bigr)
    \prec 0
\end{split}\end{equation}
as a sufficient condition for \eqref{eq:mjls_stability_vertices_lmi}.
Notably, the LMI constraints for different $i \in \mathcal{K}$ and $r \in \set{1, \ldots, n_\Gamma}$ are now decoupled (except for the shared variable $\tilde{X}$), such that it is sufficient to consider $\Gamma^{(r)}$ row by row.

The second and most important step in the process of simplifying the analysis is to assume that, in the probability distributions described by the vertices $\Gamma^{(r)}$, the processes $\{\theta_k(e^{ij})\}$ corresponding to different edges $e^{ij}$ and $e^{i'j'}$ are independent.
This implies that each row of $\Gamma^{(r)}$ is a probability vector that can be written as
\begin{equation}\label{eq:probability_map}
    f(p) \coloneqq \begin{bmatrix}
        1 - p_1 \\
        p_1
    \end{bmatrix} \otimes \begin{bmatrix}
        1 - p_2 \\
        p_2
    \end{bmatrix} \otimes \dotsb \otimes \begin{bmatrix}
        1 - p_m \\
        p_m
    \end{bmatrix}
\end{equation}
for some vector $p \in \mathcal{P}_m \coloneqq [0,1]^m$.
Now, chose two constants $0 \leq \rho_l \leq \rho_u \leq 1$ and define the set
\[\hat{\bm{t}}_m \coloneqq \set[\big]{t = f(p) \given p \in \mathcal{P}_m, p_i \in \set{\rho_l, \rho_u}} \subset \mathcal{S}_{2^m}\]
of probability vectors in which each edge has either probability $\rho_l$ or $\rho_u$ of transmitting successfully.
Combining these vectors into matrices, we obtain a set of TPMs
\[\hat{\bm{\Gamma}}_m \coloneqq \set*{\Gamma = \begin{bmatrix} t_1 & t_2 & \dotso & t_{2^m} \end{bmatrix}^\top \given t_i \in \hat{\bm{t}}_{m}} \subset \mathcal{T}_m\]
describing temporally but not spatially correlated packet loss, where $\mathcal{T}_m$ can be identified with $\mathcal{S}_{2^m} \times \dotsm \times \mathcal{S}_{2^m}$ with $2^m$ factors and is the set of all $2^m \times 2^m$ row-stochastic matrices.

Following next is the main stability condition of this paper.
It relies on treating the expectation and variance of $L_{\sigma_k}$ as an uncertainty, captured edgewise by the uncertainty set
\begin{equation*}
    \bm{\Delta} \coloneqq \set*{
        \begin{bmatrix}
            a & 0 \\
            b & 0 \\
            0 & a
        \end{bmatrix}
        \given
        \begin{aligned}
            a, b &\geq 0               \\[-0.5ex]
            a^2  &\in [\rho_l, \rho_u] \\[-0.5ex]
            b^2  &= 1 - a^2
        \end{aligned}
    }
\end{equation*}
using the full-block S-procedure \cite{Scher2000}.
The result is a variant of the approach proposed in \cite{HespeW2023} for analysing packet loss in MAS without spatial correlation.
\begin{theorem}\label{thm:bernoulli_vertices}
    The MJLS~\eqref{eq:mjls} is robustly mean-square stable for all $\Gamma \in \conv\hat{\bm{\Gamma}}_m$ if there exist $\tilde{X} \succ 0$ and
    \begin{equation}\label{eq:bernoulli_vertices_multiplier}
        P \in \set*{P = P^\top \given
            \begin{bmatrix}
                * \\ *
            \end{bmatrix}^\top \!\! P \begin{bmatrix}
                \Delta \otimes I_{n_x} \\ I_{2 n_x}
            \end{bmatrix} \succ 0 \enspace \forall \Delta \in \bm{\Delta}
        }
    \end{equation}
    such that the LMIs
    \begin{subequations}\label{eq:bernoulli_vertices_lmi}\begin{gather}
        \allowdisplaybreaks
        A^{d\top} \tilde{X} A^d - \tilde{X} \prec 0 \label{eq:bernoulli_vertices_lmi_nominal} \\
        \begin{bmatrix} * \\ * \\ * \\ \matline{1} * \\ * \\ * \end{bmatrix}^\top \!\!
        \begin{bmatrix}[ccc|c]
            -\tilde{X} & & & \\
            & \tilde{X} & & \\
            & & 2\tilde{X} & \\ \matline{4}
            & & & P
        \end{bmatrix}
        \begin{bmatrix}
            I_{n_x} & 0 \\
            A^d & \sqrt{\lambda} \mathcal{I}_1 \\
            0   & \mathcal{I}_2 \\ \matline{2}
            0   & I_{3 n_x} \\
            0   & \mathcal{I}_3 \\
            \sqrt{\lambda} A^c  & 0
        \end{bmatrix} \prec 0 \label{eq:bernoulli_vertices_lmi_modal}
    \end{gather}\end{subequations}
    hold for $\lambda \in \set{\lambda_2, \lambda_N}$, where $\lambda_2$ and $\lambda_N$ are the smallest and largest non-zero eigenvalue of $L^0$, respectively, and $\mathcal{I}_j$ is the $j$th block-row of $I_3 \otimes I_{n_x}$.
\end{theorem}
\begin{proof}
    A detailed derivation of a related result can be found in \cite{HespeW2023}.
The idea is to use the full-block S-procedure to ensure $\tilde{X}$ is a solution of \eqref{eq:mjls_stability_decomposed_lmi} for all probability vectors $t \in \hat{\bm{t}}_m$ (with $t_{ij}^{(r)}$ replaced by $t_j$).
Imposing block-repeated structure into both the Lyapunov matrix and the multiplier lets us take advantage of the decomposable systems framework to obtain constraints with the size of a single agent (see \cite{MassiV2009, HespeSDWT2023} for details).
Furthermore, due to convexity of \eqref{eq:bernoulli_vertices_lmi_modal} in $\sqrt{\lambda}$, it is sufficient to evaluate the constraint for $\lambda_2$ and $\lambda_N$.
Finally, we invoke Corollary~\ref{cor:mjls_stability_vertices} with vertices $\Gamma^{(r)} \in \hat{\bm{\Gamma}}_m$, \ie TPMs that can be assembled row-wise from $t \in \hat{\bm{t}}_m$, and conclude robust mean-square stability for the polytope $\conv \hat{\bm{\Gamma}}_m$.

\end{proof}
\begin{remark}
    In contrast to Theorem~\ref{thm:mjls_stability} and Corollary~\ref{cor:mjls_stability_vertices}, the number of variables and constraints in Theorem~\ref{thm:bernoulli_vertices} is constant, and their dimensions depend on neither $N$, $m$, nor $n_\Gamma$ but only on the state dimension $n_x$.
    The stability test can thus be applied to arbitrarily large MAS without scalability issues.
\end{remark}

\section{Characterisation of the Uncertainty Set}\label{sec:simplex}
Even though both Theorem~\ref{thm:mjls_stability} and \ref{thm:bernoulli_vertices} are testing for robust mean-square stability, the interpretation of their robustness result is different.
For Theorem~\ref{thm:mjls_stability}, the set $\bm{\Gamma}$ of uncertain TPMs is part of the system description, \ie the question of robust stability is meaningful for the MJLS only in combination with the corresponding uncertainty set.
On the other hand, Theorem~\ref{thm:bernoulli_vertices} is constructing the uncertainty set $\conv\hat{\bm{\Gamma}}_m$ in a way that is computationally favourable and up to now it is unclear which TPMs are covered.
We thus need to answer whether stabilizing the MAS for all $\Gamma \in \conv\hat{\bm{\Gamma}}_m$ provides non-trivial distributional robustness or not.

We start by studying the smallest subspace that contains $\hat{\bm{t}}_m$ and its relation to the simplex $\mathcal{S}_{2^m}$.
Before stating the result, note that there always exists a $2^m-1$ dimensional affine subspace containing $\hat{\bm{t}}_m$ since $\card{\hat{\bm{t}}_m}\, \leq\, 2^m$.
\begin{lemma}\label{lem:affine_subspace}
    Let $\rho_l < \rho_u$.
    There exists no affine subspace of dimension $2^m - 2$ that contains $\hat{\bm{t}}_m$.
\end{lemma}
\begin{proof}
    \arxivTF{Arrange the vectors in $\hat{\bm{t}}$ as
\begin{equation}\label{eq:vertex_matrix}
    M = \begin{bmatrix}
        f\bigl(p^{(1)}\bigr) & f\bigl(p^{(2)}\bigr) & \dotso & f\bigl(p^{(2^m)}\bigr)
    \end{bmatrix}
\end{equation}
where $p^{(r)} \coloneqq \rho_l \mathbf{1}_m + (\rho_u -\rho_l) \binary(r-1)$ and $\binary(\cdot)$ maps integers to their base-2 representation ordered from highest to lowest bit.
Note that $M = \tilde{M} \otimes \dotsb \otimes \tilde{M}$ with $m$ factors of
\begin{equation*}
    \tilde{M} = \begin{bmatrix}
        1 - \rho_l & 1 - \rho_u \\
        \rho_l     & \rho_u
    \end{bmatrix}
\end{equation*}
and that $\det \tilde{M} = \rho_u - \rho_l$.
By the rules for determinants of Kronecker product (\cf \cite{Steeb1991}), $\det M > 0$ follows, which implies that $M \lambda = 0$ only if $\lambda = 0$.
This guarantees that the points in $\hat{\bm{t}}_m$ are affinely independent \cite[Section 1]{Rocka1997} and that there exists no $2^m - 2$ dimensional affine subspace that contains $\hat{\bm{t}}_m$.
}{Omitted for brevity, see \cite{ExtendedVersion}.}
\end{proof}

Because of $\hat{\bm{t}}_m \subset \mathcal{S}_{2^m}$, Lemma~\ref{lem:affine_subspace} demonstrates that as long as we allow for some uncertainty in the probability, \ie $\rho_l \neq \rho_u$, the smallest affine subspace containing $\mathcal{S}_{2^m}$ is also the smallest subspace containing $\hat{\bm{t}}_m$.
Intuitively, this means the uncertainty model provided by $\hat{\bm{t}}_m$ is rich enough to describe the relevant perturbed probability vectors.

Building upon this result, we consider not just individual probabilities but intervals of probabilities.
It is thus useful to apply $f$ to whole sets, where with a slight abuse of notation we use
$f(\mathcal{M}) \coloneqq \set{f(p) \given p \in \mathcal{M}}$
for sets $\mathcal{M} \subseteq \mathcal{P}_m$.
We can then obtain a useful identity for the function~$f$:

\begin{definition}\label{def:extreme_point}
    A point $p$ in a convex set $\mathcal{M}$ is said to be an extreme point of $\mathcal{M}$ if there exist no $p', p'' \in \mathcal{M}$, $p' \neq p''$ such that $p = \beta p' + (1-\beta) p''$ with $\beta \in (0,1)$.
    We denote the set of extreme points of $\mathcal{M}$ by $\mu(\mathcal{M}) \subseteq \mathcal{M}$.
\end{definition}
\begin{lemma}\label{lem:extreme_points}
    Let $\mathcal{B} = [\rho_l, \rho_u]^m \subset \mathcal{P}_m$ with $\rho_l < \rho_u$.
    Then $f$ satisfies the identity $\mu\bigl(\conv f(\mathcal{B})\bigr) = f\bigl(\mu(\mathcal{B})\bigr)$.
\end{lemma}
\begin{proof}
    \arxivTF{Suppose for some $t \in \conv f(\mathcal{B})$ there exists no $p \in \mathcal{B}$ such that $t = f(p)$ and thus $t \in \bigl(\conv f(\mathcal{B})\bigr) \setminus f(\mathcal{B})$.
It can be shown that the extreme points of $\conv\mathcal{M}$ must be in $\mathcal{M}$ and therefore $t \notin \mu\bigl(\conv f(\mathcal{B})\bigr)$.
On the other hand, suppose there exists a $p \in \mathcal{B}$ with $t = f(p)$ but $p \notin \mu(\mathcal{B})$.
Since $\mathcal{B}$ is axis-aligned, there exist $\beta \in (0,1)$ and $p', p'' \in \mathcal{B}$ that differ in exactly one coordinate such that $p = \beta p' + (1-\beta) p''$.
From the associativity and distributivity of the Kronecker product, it then follows that
\[t = f(p) = \beta f(p') + (1-\beta) f(p''),\]
where injectivity of $f$ implies $f(p') \neq f(p'')$ and thus that $t \notin \mu\bigl(\conv f(\mathcal{B})\bigr)$ and $\mu\bigl(\conv f(\mathcal{B})\bigr) \subseteq f\bigl(\mu(\mathcal{B})\bigr)$.

For the converse, note that, due to its non-empty interior, $\card{\mu(\mathcal{B})} = 2^m$.
Furthermore, we have $\hat{\bm{t}} \subset f(\mathcal{B})$ and thus Lemma~\ref{lem:affine_subspace} implies that there exists no $2^m-2$ dimensional affine subspace that contains $f(\mathcal{B})$.
Now, suppose there exists a $p \in \mu(\mathcal{B})$ for which $f(p) \notin \mu\bigl(\conv f(\mathcal{B})\bigr)$ and thus the first half of the proof implies $\card[\big]{\mu\bigl(\conv f(\mathcal{B})\bigr)} \leq 2^m-1$.
However, every set of $2^m-1$ points (and thus their convex hull) is contained in some $2^m-2$ dimensional affine subspace, contradicting the statement above.
}{Omitted for brevity, see \cite{ExtendedVersion}.}
\end{proof}

Together, Lemma~\ref{lem:affine_subspace} and \ref{lem:extreme_points} imply that $\conv\hat{\bm{t}}_m$ is a simplex with non-empty relative interior if $\rho_l < \rho_u$, demonstrating that Theorem~\ref{thm:bernoulli_vertices} is guaranteeing stability for a non-trivial set of distributions.
The result can be strengthened by studying which probability distributions are guaranteed to be in the relative interior:

\begin{theorem}\label{thm:uncertainty_richness}
    Given $p \in \mathcal{P}_m$ and $0 \leq \rho_l' \leq \rho_l < \rho_u \leq \rho_u' \leq 1$ with corresponding sets $\hat{\bm{\Gamma}}_m$ and $\hat{\bm{\Gamma}}_m'$.
    Then,
    \begin{enumerate}[i)]
        \item $\rho_l = 0$ and $\rho_u = 1$ imply $\conv\hat{\bm{\Gamma}}_m = \mathcal{T}_m$,
        \item $\rho_l' < \rho_l$ or $ \rho_u < \rho_u'$ imply $\conv\hat{\bm{\Gamma}}_m \subset \conv\hat{\bm{\Gamma}}_m'$, and \label{itm:uncertainty_richness_monotonic}
        \item $p \in (\rho_l, \rho_u)^m$ implies $f(p) \mathbf{1}_{2^m}^\top \in \relint \conv\hat{\bm{\Gamma}}_m$. \label{itm:uncertainty_richness_interior}
    \end{enumerate}
\end{theorem}
\begin{proof}
    Because of $\conv(A \times B) = \conv(A) \times \conv(B)$ and $\relint(A \times B) = \relint(A) \times \relint(B)$ \cite{Berts2009}, we can reduce claims on $\conv\hat{\bm{\Gamma}}_m$ to $\conv\hat{\bm{t}}_m$.
Furthermore, Lemma~\ref{lem:extreme_points} with $\hat{\bm{t}}_m = f\bigl(\mu([\rho_l, \rho_u]^m)\bigr)$ imply that $\conv \hat{\bm{t}}_m = \conv f\bigl([\rho_l, \rho_u]^m\bigr)$.

\begin{enumerate}[i)]
    \item With $\rho_l = 0$ and $\rho_u = 1$, we have $\conv\hat{\bm{t}}_m = \mathcal{S}_{2^m}$.
    \item From the assumptions, it follows that $[\rho_l, \rho_u] \subset [\rho_l', \rho_u']$ and that both have non-empty interior.
        Thus Lemma~\ref{lem:extreme_points} implies $\conv f\bigl([\rho_l, \rho_u]^m\bigr) \neq \conv f\bigl([\rho_l', \rho_u']^m\bigr)$.
        We conclude $f\bigl([\rho_l, \rho_u]^m\bigr) \subset f\bigl([\rho_l', \rho_u']^m\bigr)$ from injectivity of $f$ and therefore $\conv \hat{\bm{t}}_m \subset \conv f\bigl([\rho_l', \rho_u']^m\bigr)$.
    \item We have to show $f(p) \in \relint\conv\hat{\bm{t}}_m$ and proceed by induction.
        With $m=1$, $f$ is affine in $p$.
        By assumption, there thus exists a $\beta \in (0,1)$ such that
        \[f(p) = f\bigl(\beta \rho_l + (1-\beta) \rho_u\bigr) = \beta f(\rho_l) + (1-\beta) f(\rho_u)\]
        and therefore $f(p)$ is in the relative interior.

        \begin{figure}
            \centering
            \pgfmathsetmacro{\rhol}{0.3}
\pgfmathsetmacro{\rhou}{0.85}

\pgfmathsetmacro{\pone}{0.7}
\pgfmathsetmacro{\ptwo}{0.7}

\pgfmathsetmacro{\tick}{0.05}

\tikzset{
    every node/.style={font=\scriptsize},
    point/.style={
        circle,
        fill,
        inner sep=1pt,
        outer sep=0pt,
    },
}

\NewDocumentCommand{\addcoordinate}{m >{\SplitArgument{1}{,}}r()}{\parsecoords{#1}#2}
\NewDocumentCommand{\parsecoords}{m m m}{%
    \pgfmathsetmacro{\rx}{sqrt(3)/6 * (4*#2*#3 - 1)}%
    \pgfmathsetmacro{\ry}{sqrt(6)/6 * (#2*(3 - 2*#3) - 1)}%
    \pgfmathsetmacro{\rz}{(sqrt(2)/2 * (#2 - 1)}%
    \coordinate (#1) at (\rx, \ry, \rz);
}

\subfloat[Uncertainty set in the edge probabilities.]{
    \centering
    \begin{tikzpicture}[scale=1.7, >=Stealth]
        \draw[->] (0, 0) -- (1.2, 0) node[below]{$p_1$};
        \draw[->] (0, 0) -- (0, 1.2) node[left]{$p_2$};
        
        \draw (\rhol, \tick) -- (\rhol, -\tick) node[below]{$\rho_l$};
        \draw (\rhou, \tick) -- (\rhou, -\tick) node[below]{$\rho_u$};
        \draw (\tick, \rhol) -- (-\tick, \rhol) node[left]{$\rho_l$};
        \draw (\tick, \rhou) -- (-\tick, \rhou) node[left]{$\rho_u$};
    
        \draw (\rhol, \rhol) node[point](v1){};
        \draw (\rhol, \rhou) node[point](v2){};
        \draw (\rhou, \rhou) node[point, label={above:$v^*$}](v3){};
        \draw (\rhou, \rhol) node[point, label={below:$v_1$}](v4){};
    
        \draw (v1) -- (v2);
        \draw (v3) -- (v4);
        \draw[Paired-B, thick] (v2) -- node[above, pos=0.4]{$\bar{F}$} (v3);
        \draw[Paired-B, thick] (v4) -- node[below, pos=0.6]{$F$} (v1);
        
        \draw[Paired-D, thick] (\rhol, \ptwo) -- (\rhou, \ptwo) node[right]{$F'$};
        \draw (\pone, \ptwo) node[point, label={below:$p$}](p){};
    \end{tikzpicture}
    \label{fig:probability_sketches_edges}
}
\hspace{6mm}
\subfloat[Projected view of the probability vectors mapped by $f$.]{
    \centering
    \begin{tikzpicture}[scale=3.75]
        \addcoordinate{v1}(\rhol, \rhol);
        \addcoordinate{v2}(\rhol, \rhou);
        \addcoordinate{v3}(\rhou, \rhou);
        \addcoordinate{v4}(\rhou, \rhol);
    
        \addcoordinate{p}(\pone, \ptwo);
        \addcoordinate{f1}(\rhol, \ptwo);
        \addcoordinate{f2}(\rhou, \ptwo);
    
        \draw[Paired-D, thick] (f1) -- (f2) node[above right, pos=0.93]{$f(F')$};
        \draw (p) node[point, label={[label distance=-3pt]above left:$f(p)$}]{};

        \fill[pattern=north east lines, pattern color=black!40!white] (v1) -- (v2) -- (v4) -- (v1);
        \draw ($0.5*(v1) + 0.25*(v2) + 0.25*(v4)$) node{$\mathcal{F}$};

        \draw (v1) -- (v2);
        \draw (v3) -- (v4);
        \draw[Paired-B, thick] (v2) -- node[below]{$f(\bar{F})$} (v3);
        \draw[Paired-B, thick] (v4) -- node[left, pos=0.6]{$f(F)$} (v1);
        \draw[dashed] (v1) -- (v3);
        \draw[dashed] (v2) -- (v4);
        
        \draw (v1) node[point]{};
        \draw (v2) node[point]{};
        \draw (v3) node[point, label={right:$f(v^*)$}]{};
        \draw (v4) node[point, label={above:$f(v_1)$}]{};
    \end{tikzpicture}
    \label{fig:probability_sketches_simplex}
}
            \caption{Visualization of the points and sets introduced in the proof of Theorem~\ref{thm:uncertainty_richness}, before and after being mapped by $f$.
                The four dimensional probability vectors are projected orthogonally onto $\real^3$.}
            \label{fig:probability_sketches}
        \end{figure}
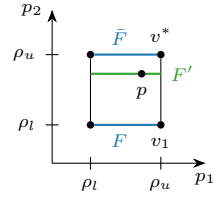
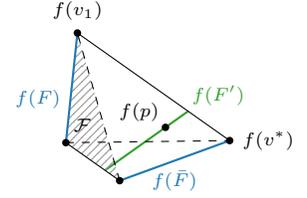

        Next, assume the statement holds for some $m \geq 1$ and consider the problem for $m+1$.
        From the hypercube $[\rho_l, \rho_u]^{m+1} \subseteq \mathcal{P}_{m+1}$, choose an arbitrary pair of opposing facets $F$ and $\bar{F}$ and in addition a vertex $v^* \in \mu(\bar{F})$.
        We denote by $V \subset \mu(F) \times \mu(\bar{F})$ the $2^m$ pairs of vertices $(v_1, v_2)$ such that $\overline{v_1 v_2}$ is an edge of $[\rho_l, \rho_u]^{m+1}$, where $\overline{v_1 v_2}$ is the line segment between $v_1$ and $v_2$.
        Finally, introduce $\mathcal{F} \coloneqq \conv \bigl(\hat{\bm{t}}_{m+1} \setminus \set{f(v^*)}\bigr)$, which by appropriate choice of $v^*$ can describe all facets of the simplex $\conv \hat{\bm{t}}_{m+1}$ that contain $f\bigl(\mu(F)\bigr)$.
        A visualization of the points and sets can be found in Fig.~\ref{fig:probability_sketches}.
        By construction, $f(v^*) \notin \mathcal{F}$ and there exists a unique $v_1 \in \mu(F)$ such that $(v_1, v^*) \in V$.
        Since $\overline{v_1 v_2}$ is axis-aligned for all $(v_1, v_2) \in V$, we thus have
        \[f\bigl(\beta v_1 + (1-\beta) v^*\bigr) = \beta f(v_1) + (1-\beta) f(v^*) \notin \mathcal{F}\]
        for all $\beta \in [0, 1)$.
        Now, consider $\beta' \in (0, 1)$ such that \[p \in F' \coloneqq \conv\set[\big]{\beta' v_1 + (1-\beta') v_2 \given (v_1, v_2) \in V},\]
        where $\beta'$ is guaranteed to exist by assumption.
        Then $p \in \relint F'$, which by the induction hypothesis implies $f(p) \in \relint \conv f(F')$.
        On the other hand, it follows from $f\bigl(\beta' v_1 + (1-\beta') v^*\bigr) \notin \mathcal{F}$ that $f\bigl(\mu(F')\bigr) \not\subset \mathcal{F}$ and therefore, because the choice of $v^*$ was arbitrary, there is no facet that contains $f\bigl(\mu(F')\bigr)$.
        We conclude that
        \[f(p) \in \relint\conv f(F') \subset \relint\conv\hat{\bm{t}}_{m+1}\]
        because $\conv f(F')$ is not entirely contained in the relative boundary of $\conv\hat{\bm{t}}_{m+1}$ \cite[Corollary~6.5.2]{Rocka1997}.
        \qedhere
\end{enumerate}

\end{proof}

Theorem~\ref{thm:uncertainty_richness} proves an important property of the uncertainty model:
For probability distributions of independent links with probability within $(\rho_l, \rho_u)$, it is guaranteed that there exists no small perturbation pushing them out of the considered uncertainty set $\conv\hat{\bm{\Gamma}}_m$.
Since the smallest affine subspaces containing $\hat{\bm{\Gamma}}_m$ and $\mathcal{T}_m$ are identical by Lemma~\ref{lem:affine_subspace}, this includes spatio-temporal correlations.
Remarkably, this is the case even though no dependent links are considered in the construction of the sets $\hat{\bm{t}}_m$ and $\hat{\bm{\Gamma}}_m$.
However, no quantification of the robustness is available, \ie there are no bounds on the least amount of robustness this result provides.

\section{Application Examples}\label{sec:example}
To demonstrate the applicability of the proposed results, this section is discussing two examples.

\subsection{Convergence of First-Order Consensus}\label{sec:example_consensus}
The first example applies Theorem~\ref{thm:bernoulli_vertices} to the first-order consensus problem (\cf \cite{OlfatFM2007}).
Therefore, we consider a set of $N$ integrators with state $x^i_{k} \in \real$ following the difference equation
\[x^i_{k+1} = x^i_k + u^i_k.\]
The agents have the shared goal of reaching agreement asymptotically and apply the consensus protocol
\[u^i_k = \kappa \sum_{j \in \mathcal{N}_i} \theta_k(e^{ij}) \big(x^j_k - x^i_k\big)\]
with gain $\kappa > 0$, where $\mathcal{N}_i \coloneqq \set{j \in \mathcal{V} \given e^{ij} \in \mathcal{E}^0}$ is the set of neighbours of agent~$i$.
Recall that we only consider undirected graphs, \ie $e^{ij} = e^{ji}$ and thus $\theta_k(e^{ij}) = \theta_k(e^{ji})$.

By stacking the states as $x_k = [x_k^1, \dotsc, x_k^N]^\top$, the consensus problem can be formulated in terms of the MJLS~\eqref{eq:mjls} as $x_{k+1} = (I_N - \kappa L_{\sigma_k}) x_k$, which can be decomposed to $A^d = 1$ and $A^c = -\kappa$.
Furthermore, introduce $\Pi \coloneqq I_N - \frac 1N \mathbf{1}_N \mathbf{1}_N^\top$ as the projection onto the disagreement space, \ie the orthogonal complement to the subspace spanned by $\mathbf{1}_N$.
This leads to the following guarantee on the convergence of the agents:

\begin{proposition}\label{pro:consensus}
    Let $\mathcal{G}^0$ be connected, $0 < \kappa < \frac{2}{\lambda_N}$, and $0 < \rho_l \leq \rho_u = 1$, where $\lambda_N$ is the largest eigenvalue of $L^0$.
    Then the first-order agents applying the consensus protocol satisfy
    \[\lim_{k \to \infty} \expect[\big]{\norm{\Pi x_k}^2} = 0\]
    for all $\Gamma \in \hat{\bm{\Gamma}}_m$ and all initial conditions $x_0$, and $\sigma_0$.
\end{proposition}
\begin{proof}
	\arxivTF{Since $\mathcal{G}^0$ is undirected, there exists an orthogonal transformation $U^\top U = I_N$ such that $U^\top L^0 U$ is diagonal, and furthermore $U$ can be chosen as $U = [\mathbf{1}_N / \sqrt{N} \enspace \tilde{U}]$ such that $\Pi U = [0 \enspace \tilde{U}]$.
Applying the change of variables $\tilde{x}_k \coloneqq U^\top x_k$, we obtain
\begin{equation*}
    \tilde{x}_{k+1} = \begin{bmatrix}
        A^d & 0 \\
        0   & I_{N-1} \otimes A^d + \tilde{L}_{\sigma_k} \otimes A^c
    \end{bmatrix} \tilde{x}_k
\end{equation*}
with $\tilde{L}_i \coloneqq \tilde{U}^\top L_i \tilde{U}$.
From $\Pi x_k = [0 \enspace \tilde{U}] \tilde{x}_k$ it follows that the first entry of $\tilde{x}_k$ can neglected, which amounts to disregarding \eqref{eq:bernoulli_vertices_lmi_nominal} in Theorem~\ref{thm:bernoulli_vertices}.
We thus only have to prove the existence of $\tilde{X} \succ 0$ and $P$ satisfying \eqref{eq:bernoulli_vertices_multiplier} and \eqref{eq:bernoulli_vertices_lmi_modal}. 

We proceed by constructing a solution to the LMIs that is valid for all $\rho_l > 0$.
First choose an $0 < \varepsilon < \frac{2 - \kappa \lambda_N}{4}$, which exists by assumption.
Then, set $\tilde{X} = \kappa$ and
\begin{align*}
    P = \begin{bmatrix}
        Q & S^\top \\
        S & R
    \end{bmatrix} & & \text{with} & &
    S = \begin{bmatrix}
        0 & 0 & -1 \\
        1 & 0 &  0
    \end{bmatrix}
\end{align*}
and where $Q, R$ are diagonal matrices with entries
\begin{align*}
    Q_{11} &= Q_{22} = -\kappa \lambda_N - \varepsilon, & R_{11} &= \kappa \lambda_N + 2\varepsilon, \\
    Q_{33} &= 2 - \kappa \lambda_N - 3\varepsilon,       & R_{22} &= \rho_L(\kappa \lambda_N - 2 + 4\varepsilon).
\end{align*}
Condition~\eqref{eq:bernoulli_vertices_multiplier} simplifies to
$R_{11} + \rho Q_{11} + (1-\rho) Q_{22} > 0$ and
$R_{22} + \rho Q_{33} > 0$ for $\rho \in \set{\rho_l, 1}$, while \eqref{eq:bernoulli_vertices_lmi_modal} results in $R_{22} < 0$, $Q_{11} + \kappa \lambda < 0$, $Q_{22} + 2\kappa < 0$ and $Q_{33} + R_{11} < 2$ for $\lambda \in \set{\lambda_2, \lambda_N}$.
Since $\lambda_N \geq 2$ for connected graphs \cite{GroneM1994}, the proposed $P$ satisfies all these conditions.
}{%
		The result is shown by explicitly construction a solution $X \succ 0$ and $P$ to \eqref{eq:bernoulli_vertices_multiplier} and \eqref{eq:bernoulli_vertices_lmi_modal} that is valid for all $\rho_l > 0$.
		A detailed construction can be found in \cite{ExtendedVersion}.
	}%
\end{proof}

Proposition~\ref{pro:consensus} establishes that applying the consensus protocol leads to the agents asymptotically reaching agreement under the same conditions as with ideal communication \cite{OlfatFM2007} if every link has a non-zero chance to transmit.
Similar statements have been reported in \cite[Chapter 5]{MesbaE2010} and \cite{ZhangT2012}, relying however on the independence of the transmissions.
This demonstrates that Theorem~\ref{thm:bernoulli_vertices} can be used to obtain stability guarantees for broad ranges of probability distributions.

\subsection{Fitting Simulated Traffic Communication Data}\label{sec:example_traffic}
The second example is concerned with demonstrating that the uncertainty set $\conv\hat{\bm{\Gamma}}_m$ not only has the abstract properties guaranteed by Theorem~\ref{thm:uncertainty_richness} but covers TPMs encountered in practical applications.
As our dataset, we use the TPMs obtained in \cite{RazzaDSZWFVF2022} through means of high-fidelity simulations of vehicle-to-vehicle communication.
We can then fit our model to the data by utilizing that, given a TPM~$\Gamma$ and an interval $[\rho_l, \rho_u]$, $\Gamma \in \conv\hat{\bm{\Gamma}}_m$ is equivalent to the existence of a solution $\Lambda \in \real^{2^m \times 2^m}$ to the linear feasibility problem
\begin{subequations}\label{eq:tpm_lp}\begin{align}\SwapAboveDisplaySkip
    \Lambda M^\top            &= \Gamma \\
    \Lambda \mathbf{1}_{2^m}  &= \mathbf{1}_{2^m} \\
    \Lambda                   &\geq 0, \label{eq:tpm_lp_positivity}
\end{align}\end{subequations}
where $M$ is defined \arxivTF{in \eqref{eq:vertex_matrix} }{as
	\[M \coloneqq \begin{bmatrix}f\bigl(p^{(1)}\bigr) & f\bigl(p^{(2)}\bigr) & \dotso & f\bigl(p^{(2^m)}\bigr)\end{bmatrix},\]
	$f\bigl(p^{(i)}\bigr)$ are the elements of $\hat{\bm{t}}_m$,
}%
and \eqref{eq:tpm_lp_positivity} is an element-wise inequality.
We can thus obtain the smallest interval $[\rho_l, \rho_u]$ such that $\Gamma \in \conv\hat{\bm{\Gamma}}_m$ through iterative solutions of \eqref{eq:tpm_lp}, which we implemented using \textsc{Yalmip} \cite{Loefb2004}.

The dataset of \cite{RazzaDSZWFVF2022} is parametrized through two quantities, the traffic density and the inter-packet gap (IPG), \ie the time between two received packets.
The results of the interval minimization are listed in Table~\ref{tab:intervals}.
\begin{table}
    \centering
    \caption{Smallest probability intervals such that the TPM is contained in $\hat{\bm{\Gamma}}_m$ for six platooning scenarios from \cite{RazzaDSZWFVF2022}}
    \label{tab:intervals}

	\pgfplotstableread[col sep=comma]{data/minimal_interval.csv} \datatable

	\pgfplotstabletypeset[
	    create on use/diff/.style={
            create col/expr={\thisrow{upper}-\thisrow{lower}}
        },
	    columns={density, ipg, lower, upper, diff},
        every column/.style={
			column type=r,
			fixed,
			zerofill,
		},
		columns/density/.style={
			string type,
			column name={Traffic Density},
			column type=l,
		},
		columns/ipg/.style={
		    string type,
			column name={IPG in \si{\milli\second}},
			string replace*={-}{ -- },
		},
	    columns/lower/.style={column name=$\rho_l$},
		columns/upper/.style={column name=$\rho_u$},
		columns/diff/.style={column name=$\Delta \rho$},
        every head row/.style={
			before row={
                \toprule
                \multicolumn{2}{l}{Scenario} & \multicolumn{3}{r}{Probability Interval} \\ \cmidrule(r){1-2} \cmidrule(l){3-5}
            },
			after row=\midrule
		},
		every last row/.style={after row=\bottomrule},
	] \datatable
\end{table}
The table shows that we are able to find non-trivial intervals for all six cases, which demonstrates that it is possible to capture TPMs of correlated networks as convex combinations of probability distributions consisting of independent links.

\section{Conclusions}\label{sec:conclusions}
This paper has presented a scalable condition for robust mean-square stability of networked MASs subject to stochastic packet loss with spatio-temporal correlations.
Its main novelty is that the stability condition is derived from distributions of independent communication links only, leading to conditions with a computational complexity that is independent of the number of agents.
Nonetheless, it was shown that the stabilized uncertainty set has non-empty relative interior, providing distributional robustness.
An important question left to study is how to obtain a quantitative measure for the robustness provided by the result, that is, an explicit lower bound on the available robustness margin.

\bibliographystyle{IEEETran}
\bibliography{arxiv}

\end{document}